\documentclass[letterpaper, 10 pt, conference]{ieeeconf}%
\usepackage{graphics}
\usepackage{epsfig}
\usepackage{cite}
\usepackage{breakurl}
\usepackage{amsmath}
\usepackage{bm}
\usepackage{amssymb}
\usepackage{amsfonts}
\usepackage{graphicx}
\usepackage{makecell,multirow,diagbox}
\usepackage{footnote}
\usepackage[normalem]{ulem}
\usepackage{color}
\usepackage{dsfont}
\usepackage[hyphens]{url}
\usepackage[hidelinks]{hyperref}
\usepackage[ruled,linesnumbered]{algorithm2e}
\usepackage[dvipsnames]{xcolor}%
\providecommand{\U}[1]{\protect\rule{.1in}{.1in}}
\IEEEoverridecommandlockouts
\newtheorem{innercustomdef}{Definition}
\newenvironment{customdef}[1]
  {\renewcommand\theinnercustomdef{#1}\innercustomdef}
  {\endinnercustomdef}

\newtheorem{innercustomcor}{Corollary}

\newtheorem{innercustomrem}{Remark}
\newenvironment{customrem}[1]
  {\renewcommand\theinnercustomrem{#1}\innercustomrem}
  {\endinnercustomcor}

\allowdisplaybreaks[0]
\hypersetup{breaklinks=true}
\usepackage[skip=0.5pt]{caption}
\begin{document}

\title{{\LARGE \textbf{A Decentralized Optimal Control Framework for Connected
Automated Vehicles at Urban Intersections with Dynamic Resequencing}}}
\author{Yue Zhang, and Christos G. Cassandras \thanks{Supported in part by NSF under grants ECCS-1509084, CNS-1645681, and
IIP-1430145, by AFOSR under grant FA9550-15-1-0471, by DOE under grant
DOE-46100, by MathWorks and by Bosch.} \thanks{Y. Zhang and C.G. Cassandras
are with the Division of Systems Engineering and Center for Information and
Systems Engineering, Boston University, Boston, MA 02215 USA (e-mail:
joycez@bu.edu; cgc@bu.edu).} }
\maketitle


\begin{abstract}
Earlier work has established a decentralized framework to optimally control
Connected Automated Vehicles (CAVs) crossing an urban intersection without
using explicit traffic signaling while following a strict First-In-First-Out
(FIFO) queueing structure. The proposed solution minimizes energy consumption
subject to a FIFO-based throughput maximization requirement. In this paper, we
extend the solution to account for asymmetric intersections by
relaxing the FIFO constraint and including a dynamic resequencing process so
as to maximize traffic throughput. To investigate the tradeoff between
throughput maximization and energy minimization objectives, we exploit several
alternative problem formulations. In addition, the computational complexity of
the resequencing process is analyzed and proved to be bounded, which makes the
online implementation computationally feasible. The effectiveness of the
dynamic resequencing process in terms of throughput maximization is
illustrated through simulation.

\end{abstract}

\thispagestyle{empty} \pagestyle{empty}


\section{Introduction}

\label{intro}

To date, traffic light control is the prevailing method for controlling the
traffic flow in an urban area. Recent technological developments (e.g.,
\cite{Fleck2015}) have exploited data-driven control and optimization
approaches and enabled the adaptive control of traffic light cycles, which
reduces the travel delay. However, in addition to the obvious infrastructure
cost,
safety issues, e.g., rear-end collisions, often arise under traffic light
control. These issues have motivated research efforts to explore new
approaches capable of enabling a smoother traffic flow while also improving safety.

Connected and Automated Vehicles (CAVs) have the potential to drastically
improve a transportation network's performance by assisting drivers in making
better decisions, ultimately reducing energy consumption, air pollution,
congestion and accidents. One of the very early efforts exploiting the benefit
of CAVs was proposed in \cite{Levine1966}, where an optimal linear feedback
regulator is introduced for the merging problem to control a single string of
vehicles. More recently, several research efforts have been reported in the
literature for CAV coordination at intersections. Dresner and Stone
\cite{Dresner2004} proposed a reservation-based scheme for centralized
automated vehicle intersection management. Since then, numerous research
efforts have explored safe and efficient control strategies, e.g.,
\cite{Dresner2008, DeLaFortelle2010, Huang2012}. Some approaches have focused
on coordinating vehicles so as to reduce travel delay and increase
intersection throughput, e.g., \cite{Li2006,Yan2009,Zhu2015} and some have
studied intersections as polling systems \cite{Miculescu2014} so as to
determine a sequence of times assigned to vehicles on each road. Reducing
energy consumption is one of the desired objectives which has been considered in recent literature \cite{gilbert1976vehicle,
hooker1988optimal, hellstrom2010design, li2012minimum}. A detailed discussion
of the overall research in this area can be found in \cite{Rios-Torres}.

Our earlier work \cite{ZhangMalikopoulosCassandras2016} has established a
decentralized optimal control framework for coordinating on line a continuous
flow of CAVs crossing an urban intersection without using explicit traffic
signaling. For each CAV, an energy minimization optimal control problem is
formulated where the time to cross the intersection is first determined
through a throughput maximization problem. We also established conditions
under which feasible solutions to the optimal control problem exist.

The crossing sequence for the CAVs based on which the throughput maximization
problem in \cite{Malikopoulos2016} is formulated adopts a strict
First-In-First-Out (FIFO) queueing structure. This can be effective when the
intersection is physically symmetrical and the vehicle arrival rates at all
intersection entries do not differ much. However, when the intersection is
asymmetrical, the FIFO queueing structure may lead to poor scheduling and
possible congestion. Even with a fully symmetrical intersection, a strict FIFO
crossing sequence is conservative in the sense that it prevents the
intersection from further exploiting the benefits of CAVs and achieving
traffic throughput maximization.
Hence, it is necessary to design a coordination algorithm for CAVs to maximize
the traffic throughput. Zohdy et al. \cite{Zohdy2012} presented an approach
based on Cooperative Adaptive Cruise Control (CACC) for minimizing
intersection delay and hence maximizing the throughput. Lee and Park
\cite{Lee2012} considered minimizing the overlap between vehicle positions. In
this paper, we extend the optimal control solution in \cite{Malikopoulos2016}
by relaxing the FIFO constraint and including a dynamic resequencing process
so as to maximize traffic throughput. 

The paper is structured as follows. In Section II, we review the model in
\cite{ZhangMalikopoulosCassandras2016} and its generalization in
\cite{Malikopoulos2016}. In Section III, we extend the solution in
\cite{Malikopoulos2016} by relaxing the FIFO constraint and including a
dynamic resequencing process so as to maximize traffic throughput. We consider
several alternative problem formulations in order to investigate the tradeoff
between throughput maximization and energy minimization objectives. In Section
IV, we analyze the computational complexity of the resequencing process and
show it to be bounded and, on average, limited to the number of lanes in the
intersection. Conclusions and future work are given in Section V.

\section{The Model}

\label{sec:1a}

\begin{figure}[ptb]
\centering
\includegraphics[width=0.7\columnwidth]{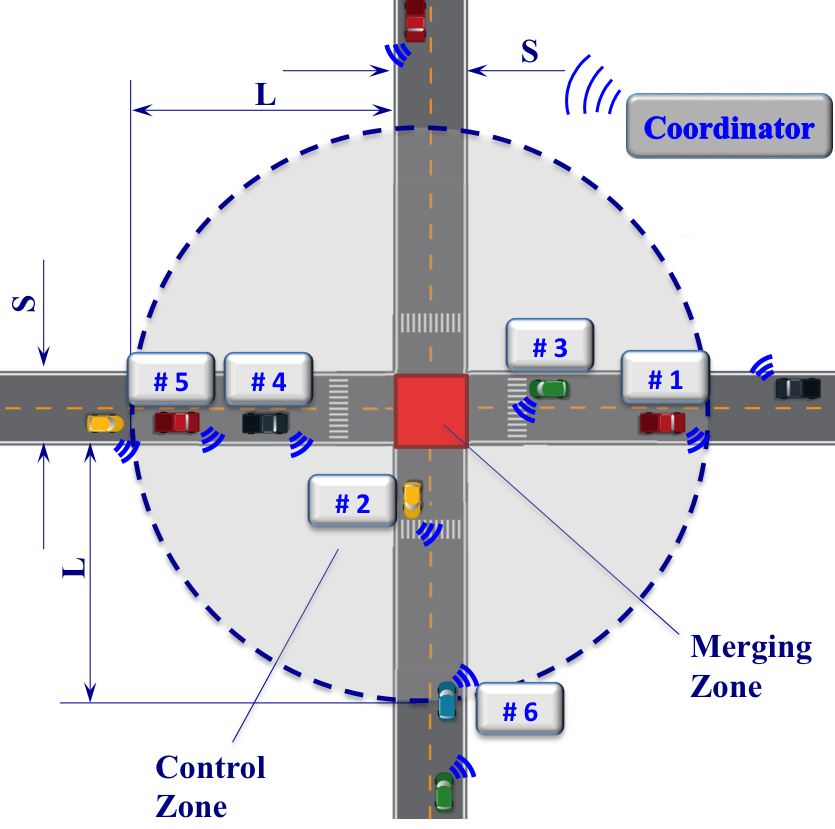} \caption{Connected
Automated Vehicles crossing an urban intersection.}%
\label{fig:intersection}%
\end{figure}

The model introduced in \cite{ZhangMalikopoulosCassandras2016} and
\cite{Malikopoulos2016} is briefly reviewed. We consider an intersection (Fig.
\ref{fig:intersection}) where the region at the center of each intersection,
called \emph{Merging Zone} (MZ), is the area of potential lateral CAV
collision and assumed to be a square of side $S$. The intersection has a
\emph{Control Zone} (CZ) and a coordinator that can communicate with the CAVs
traveling within it. The road segment from the CZ entry to the CZ exit (i.e.,
the MZ entry) is referred as a CZ segment. The length of CZ segment is $L>S$,
and it is assumed to be the same for all entry points to a given CZ.

Let $N(t)\in\mathbb{N}$ be the cumulative number of CAVs which have entered
the CZ by time $t$ and formed a queue that designates the crossing sequence in
which these CAVs will enter the MZ. There is a number of ways to form the
queue. In \cite{ZhangMalikopoulosCassandras2016} and \cite{Malikopoulos2016},
a strict First-In-First-Out (FIFO) crossing sequence is assumed, that is, when
a CAV reaches the CZ, the coordinator assigns it an integer value $i=N(t)+1$.
If two or more CAVs enter a CZ at the same time, then the corresponding
coordinator selects randomly the first one to be assigned the value $N(t)+1$.

For simplicity, we assume that each CAV is governed by a second order dynamics%
\begin{equation}
\dot{p}_{i}=v_{i}(t)\text{, }~p_{i}(t_{i}^{0})=0\text{; }~\dot{v}_{i}%
=u_{i}(t)\text{, }v_{i}(t_{i}^{0})\text{ given} \label{eq:model2}%
\end{equation}
where $p_{i}(t)\in\mathcal{P}_{i}$, $v_{i}(t)\in\mathcal{V}_{i}$, and
$u_{i}(t)\in\mathcal{U}_{i}$ denote the position, i.e., travel distance since
the entry of the CZ, speed and acceleration/deceleration (control input) of
each CAV $i$. The sets $\mathcal{P}_{i}$, $\mathcal{V}_{i}$ and $\mathcal{U}%
_{i}$ are complete and totally bounded sets of $\mathbb{R}$. These dynamics
are in force over an interval $[t_{i}^{0},t_{i}^{f}]$, where $t_{i}^{0}$ and
$t_{i}^{f}$ are the times that the vehicle $i$ enters the CZ and exits the MZ respectively.

To ensure that the control input and vehicle speed are within a given
admissible range, the following constraints are imposed:
\begin{equation}%
\begin{split}
u_{i,min}  &  \leq u_{i}(t) \leq u_{i,max},\quad\text{and}\\
0  &  \leq v_{min} \leq v_{i}(t) \leq v_{max},\quad\forall t\in\lbrack
t_{i}^{0},t_{i}^{f}].
\end{split}
\label{speed_accel constraints}%
\end{equation}

\begin{customdef}{1}
Depending on its physical location inside the CZ, CAV $i-1\in\mathcal{N}(t)$
belongs to only one of the following four subsets of $\mathcal{N}(t)$ with
respect to CAV $i$: 1) $\mathcal{R}_{i}(t)$ contains all CAVs traveling on the
same road as $i$ and towards the same direction but on different lanes, 2)
$\mathcal{L}_{i}(t)$ contains all CAVs traveling on the same road and lane as
vehicle $i$ (e.g., $\mathcal{L}_{5}(t)$ contains CAV \#4 in Fig.
\ref{fig:intersection}), 3) $\mathcal{C}_{i}(t)$ contains all CAVs traveling
on different roads from $i$ and having destinations that can cause collision
at the MZ (e.g., $\mathcal{C}_{6}(t)$ contains CAV \#5 in Fig.
\ref{fig:intersection}), and 4) $\mathcal{O}_{i}(t)$ contains all CAVs
traveling on the same road as $i$ and opposite destinations that cannot,
however, cause collision at the MZ (e.g., $\mathcal{O}_{4}(t)$ contains CAV
\#3 in Fig. \ref{fig:intersection}). \label{def:2}
\end{customdef}

To ensure the absence of any rear-end collision throughout the CZ, we impose
the \emph{rear-end safety constraint}:
\begin{equation}
s_{i}(t)=p_{k}(t)-p_{i}(t)\geq\delta,~\forall t\in\lbrack t_{i}^{0},t_{i}%
^{f}],~k\in\mathcal{L}_{i}(t) \label{rearend}%
\end{equation}
where $k$ is the CAV physically ahead of $i$ on the same lane, $s_{i}(t)$ is
the inter-vehicle distance between $i$ and $k$, and $\delta$ is the
\emph{minimal safety following distance} allowable.

A lateral collision involving CAV $i$ may occur only if some CAV $j \neq i$
belongs to $\mathcal{C}_{i}(t)$. This leads to the following definition:

\begin{customdef}{2}
For each CAV $i\in\mathcal{N}(t)$, we define the set $\Gamma_{i}$ that
includes all time instants when a lateral collision involving CAV $i$ is
possible: $\Gamma_{i}\triangleq\Big\{t~|~t\in\lbrack t_{i}^{m},t_{i}%
^{f}]\Big\}$, where $t_{i}^{m}$ is the time that CAV $i$ enters the MZ.
Consequently, to avoid a lateral collision for any two vehicles $i,j\in
\mathcal{N}(t)$ on different roads, the following constraint should hold
\begin{equation}
\Gamma_{i}\cap\Gamma_{j}=\varnothing,\text{ \ \ \ }\forall t\in\lbrack
t_{i}^{m},t_{i}^{f}]\text{, \ }j\in\mathcal{C}_{i}(t). \label{eq:lateral}%
\end{equation}

\end{customdef}

As part of safety considerations, we impose the following assumptions:
For CAV $i$, none of the constraints (\ref{speed_accel constraints}%
)-(\ref{rearend}) is active at $t_{i}^{0}$. 
The speed of the CAVs inside the MZ is constant, i.e., $v_{i}(t)=v_{i}^{m}$,
$\ \forall t\in\lbrack t_{i}^{m},t_{i}^{f}]$. This implies that $t_{i}%
^{f}=t_{i}^{m}+\frac{S}{v_{i}^{m}}$. 
Each CAV $i$ has proximity sensors and can measure local information without
errors or delays. 

The objective of each CAV is to derive an optimal acceleration/deceleration
profile, in terms of minimizing energy consumption, inside the CZ while avoiding congestion between
the two intersections. Since the coordinator is not involved in any decision
making process on the vehicle control, we can formulate $N(t)$ decentralized
tractable problems that can be solved online, that is,
\begin{gather}
\min_{u_{i}\in U_{i}}\frac{1}{2}\int_{t_{i}^{0}}^{t_{i}^{m}}K_{i}\cdot
u_{i}^{2}~dt\nonumber\\
\text{subject to}:\eqref{eq:model2},(\ref{speed_accel constraints}),t_{i}%
^{m},\text{ }p_{i}(t_{i}^{0})=0\text{, }\label{eq:decentral}\\ p_{i}(t_{i}^{m})=L,
\text{and given }t_{i}^{0}\text{, }v_{i}(t_{i}^{0}),\nonumber
\end{gather}
where $K_{i}$ is a factor to capture CAV diversity ($K_{i}=1$ for
simplicity). Note that the terminal speed $v_{i}^{m}$ is undefined and obtained from the energy minimization problem.

The terminal times for CAVs entering the MZ, i.e., $t_{i}^{m}$, can be
obtained as the solution to a throughput maximization problem based on a FIFO
crossing sequence subject to rear-end and lateral collision avoidance
constraints inside the MZ. As shown in \cite{Malikopoulos2016}, the terminal
time of CAV $i$ (i.e., $t_{i}^{m}$) can be recursively determined through%
\begin{equation}
t_{i}^{m^{\ast}}=\left\{
\begin{array}
[c]{ll}%
t_{1}^{m^{\ast}} & \mbox{if $i=1$}\\
\text{max }\{t_{i-1}^{m^{\ast}},t_{k}^{m^{\ast}}+\frac{\delta}{v_{k}^{m}%
},t_{i}^{c}\} & \text{if }i-1\in\mathcal{R}_{i}\cup\mathcal{O}_{i}\\
\text{max }\{t_{i-1}^{m^{\ast}}+\frac{\delta}{v_{i-1}^{m}},t_{i}^{c}\} &
\mbox{if $i-1\in\mathcal{L}_{i}$}\\
\text{max }\{t_{i-1}^{m^{\ast}}+\frac{S}{v_{i-1}^{m}},t_{i}^{c}\} &
\mbox{if $i-1\in\mathcal{C}_{i}$}
\end{array}
\right.  \label{def:tf}%
\end{equation}
where $t_{i}^{c}=t_{i}^{1}\mathds{1}_{v_{i}^{m}=v_{max}}+t_{i}^{2}%
(1-\mathds{1}_{v_{i}^{m}=v_{max}})$ and $
t_{i}^{1}   =t_{i}^{0}+\frac{L}{v_{max}}+\frac{(v_{max}-v_{i}^{0})^{2}%
}{2u_{i,max}v_{max}}$, $t_{i}^{2}   =t_{i}^{0}+\frac{[2Lu_{i,max}+(v_{i}^{0})^{2}]^{1/2}-v_{i}^{0}}{u_{i,max}}$.
Here, $t_{i}^{c}$ is a lower bound of $t_{i}^{m}$ regardless of the solution
of the throughput maximization problem.

An analytical solution of problem \eqref{eq:decentral} may be obtained through
a Hamiltonian analysis. Assuming that all constraints are satisfied upon
entering the CZ and that they remain inactive throughout $[t_{i}^{0},t_{i}%
^{m}]$, the optimal control input (acceleration/deceleration) over
$t\in\lbrack t_{i}^{0},t_{i}^{m}]$ is given by
\begin{equation}
u_{i}^{\ast}(t)=a_{i}t+b_{i} \label{eq:20}%
\end{equation}
where $a_{i}$ and $b_{i}$ are constants of integration. Using (\ref{eq:20}) in
the CAV dynamics \eqref{eq:model2}, the optimal speed and position are obtained as
\begin{equation}
v_{i}^{\ast}(t)=\frac{1}{2}a_{i}t^{2}+b_{i}t+c_{i} \label{eq:21}%
\end{equation}%
\begin{equation}
p_{i}^{\ast}(t)=\frac{1}{6}a_{i}t^{3}+\frac{1}{2}b_{i}t^{2}+c_{i}t+d_{i},
\label{eq:22}%
\end{equation}
where $c_{i}$ and $d_{i}$ are constants of integration. The coefficients
$a_{i}$, $b_{i}$, $c_{i}$, $d_{i}$ can be obtained given initial and terminal conditions.

Note that the analytical solution \eqref{eq:20} is valid while none of the
constraints becomes active for $t \in[t_{i}^{0}, t_{i}^{m}]$. Otherwise, the
optimal solution should be modified considering the active constraints as
discussed in \cite{Malikopoulos2016}. Also note that this formulation
\eqref{eq:decentral} does not include the safety constraint (\ref{rearend}).
The conditions under which the rear-end collision avoidance constraint does
not become active inside the CZ are provided in \cite{Zhang2016}, where it is
also shown how they can be enforced through an appropriately designed
\textit{Feasibility Enforcement Zone} that precedes the CZ.

\section{Dynamic Resequencing of Connected Automated Vehicles}

\begin{figure}[tbh]
\centering
\includegraphics[width=1\columnwidth]{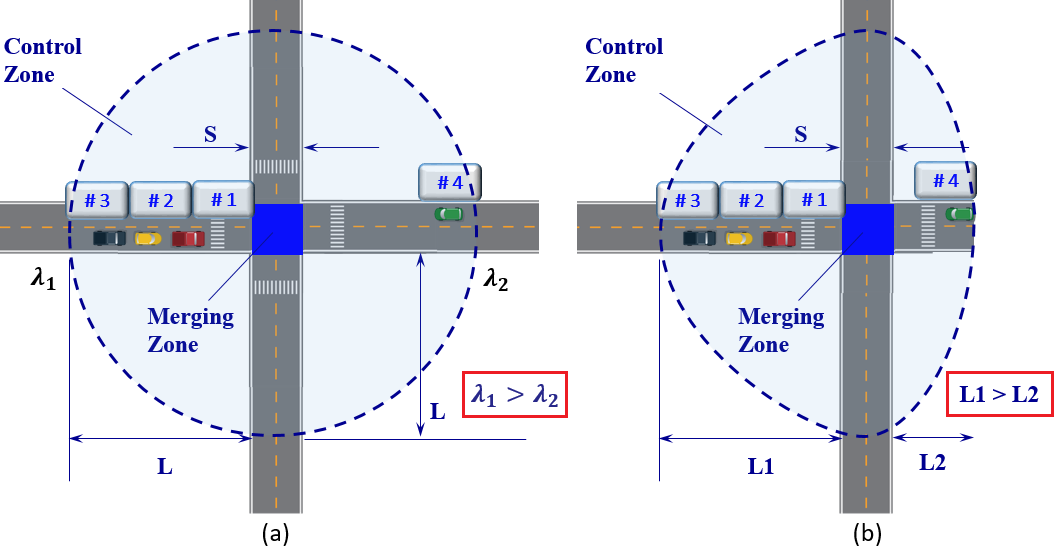} \caption{Connected
Automated Vehicles crossing an asymmetrical urban intersection.}%
\label{asym}%
\end{figure}

The crossing sequence for the CAVs based on which the throughput maximization problem  in  \cite{Malikopoulos2016} is
formulated adopts a strict FIFO queueing structure. This can be effective
when the CZ is physically symmetrical and the vehicle arrival rates at all CZ entries do not differ
much. However, when the CZ is asymmetrical (see Fig. \ref{asym}), the FIFO
queueing structure may lead to poor scheduling and possible congestion. For
example, in Fig. \ref{asym}(a) where the CZ is asymmetrical in terms of the
vehicle arrival rates, CAV \#4 entering the intersection from a CZ segment
with a lower arrival rate should wait under FIFO for the first three CAVs
crossing the MZ, which leads to unnecessary travel delay and extra energy
consumption. In Fig. \ref{asym}(b) where the CZ is asymmetrical in terms of
the physical lengths of CZ segments, CAV \#4 enters the intersection from a
shorter CZ segment and it is closer to the MZ entry, while \#4 has to
decelerate in order to let the CAV \#1, \#2 and \#3 cross the MZ first. This
again will increase travel delay. Even with a fully symmetrical CZ, a strict
FIFO queueing structure is conservative in the sense that it prevents the CZ
from achieving higher traffic throughput. For example, a CAV with higher
initial speed may tend to cross the MZ before another CAV which arrives at the
CZ earlier but with lower initial speed.

\subsection{Feasible Crossing Sequence}

A natural approach dealing with the sequencing issue is to dynamically
resequence the CAVs when a new one enters the CZ. The resequencing policy can
be position-based, i.e., the CAV closer to the MZ entry is prioritized to cross
it. Alternatively, the crossing sequence can be determined based on the
estimated travel time to the MZ. However, these methods may not be fair since
CAVs entering from the shorter CZ segment are always prioritized over those
entering from the longer CZ segment, which leads to congestion on the longer
CZ segment. A better approach is to evaluate all feasible crossing sequences
whenever a new CAV enters the CZ and select the one that maximizes traffic throughput.

Thus, our objective is to assign each arriving CAV an appropriate order to
maximize traffic throughput while maintaining the relative order of the
remaining CAVs. The problem reduces to finding all feasible crossing
sequences, computing the corresponding terminal times recursively as in
\eqref{def:tf}, and determining the one providing maximal throughput.

The first step is to find all the feasible crossing sequences, which is
equivalent to finding all the feasible orders which can be assigned to CAV
$i=N(t)+1$. Recalling that CAV $k\in\mathcal{L}_{i}$ is the vehicle physically
ahead of $i$ on the same lane ($k=0$ if such a CAV does not exist), we define
a function $f(m,n)$ which swaps the order of CAVs $m$ and $n$, i.e., after
$f(m,n)$ is evaluated, CAVs $m$ and $n$ become CAVs $n$ and $m$ respectively.
Denoting a crossing sequence as $s$ and the set containing all the feasible
crossing sequences as $\mathcal{S}_{i}$ when CAV $i$ enters the CZ, the
algorithm for deriving all feasible crossing sequences $\mathcal{S}_{i}$ is
presented as follows.

\begin{algorithm}
set $j_i := i$\;
\While{$j_i \neq k$}{
obtain a new $s$\;
\eIf{$s$ is feasible}{
add $s$ to $\mathcal{S}_i$\;
}{
break\;
}
execute $f(j_i, j_i-1)$ $\Rightarrow j_i := j_i-1$\;
}
\caption{Find the feasible crossing sequence set $\mathcal{S}_i$}
\label{feasible_sequence}
\end{algorithm}

Note that CAV $i$ cannot overtake the preceding CAV $k$ on the same lane.
Therefore, the algorithm will stop if faced with an order swap of CAVs $i$ and
$k$. After each call of $f$, the original order of $i$ which is $o(i)=i$ is
assigned a new order $o^{\prime}(i)=$ $j_{i}$, where the subscript $i$
represents the original order, and the coordinator will obtain a new crossing
sequence $s$. Recalling that $t_{i}^{c}$ is the lower bound of $t_{i}^{m}$, the
crossing sequence $s$ can only be feasible if $t_{j_{i}}^{m}\geq t_{i}^{c}$ is
satisfied, where $t_{j_{i}}^{m}$ is given through \eqref{def:tf}. Clearly, all
existing CAVs whose order is affected by the resequencing process may only
arrive at the MZ later than the original terminal times; therefore, it is not
possible for them to violate the lower bound. If $t_{j_{i}}%
^{m}<t_{i}^{c}$ holds, the sequence $s$ is infeasible and the
algorithm can terminate. This is because for any order $j_{i}^{^{\prime}%
}<j_{i}$, we have $t_{j_{i}^{^{\prime}}}^{m} \leq t_{j_{i}}^{m}<t_{i}^{c}$, hence,
there is no need to continue the algorithm. If the crossing sequence $s$ is
feasible, the coordinator will record this sequence and add it to the set
$\mathcal{S}_{i}$. This process repeats until $f$ can no longer be executed.
Note that the set $\mathcal{S}_{i}$ must be non-empty since $j_{i}=i$ itself
is always a feasible order for CAV $i$. The computational complexity of this
process will be discussed in Section IV.

\subsection{Throughput Maximization Problem Formulation}

For each feasible crossing sequence $s$ in $\mathcal{S}_{i}$, we can determine
the terminal time for each CAV iteratively through \eqref{def:tf} and obtain a
terminal time sequence $\mathbf{t}_{(2:i)}=[t_{2}^{m},\cdots,t_{i}^{m}]$. As
in \cite{ZhangMalikopoulosCassandras2016} and \cite{Malikopoulos2016}, we aim at minimizing the gaps between the terminal
times of two adjacent CAVs $i$ and $i-1$ in the sequence. Given the recursive
structure of the terminal times, this objective is equivalent to minimizing
$t_{i}^{m}-t_{1}^{m}$. Thus, our objective is
\begin{gather}
\min_{s\in\mathcal{S}_{i}}\sum_{j=2}^{i}\Big(t_{j}^{m}-t_{j-1}^{m}%
\Big)=\min_{s\in\mathcal{S}_{i}}\Big(t_{i}^{m}-t_{1}^{m}\Big)\\
\text{subject to:~}%
\eqref{eq:model2},\eqref{speed_accel constraints}, \eqref{eq:lateral},\nonumber\\
s_{i}(t)=p_{k}(t)-p_{i}(t)\geqslant\delta,~\forall t\in\lbrack t_{i}^{m}%
,t_{i}^{f}]\text{, \ }k\in\mathcal{L}_{i}(t).\nonumber\label{tt1}%
\end{gather}

Observe that $t_{1}^{m}$ is not included in the terminal time sequence since
its selection is subject to a degree of freedom reflecting the tradeoff
between energy minimization and throughput maximization. In our earlier work
\cite{ZhangMalikopoulosCassandras2016} and \cite{Malikopoulos2016}, CAV \#1 is
assumed to cruise at its initial speed so that $t_{1}^{m}=t_{1}^{0}+\frac{L}{v_{1}^{0}}$ and its
terminal speed is $v_{1}^{m}=v_{1}^{0}$. However, with resequencing, several
alternatives are possible as discussed in the sequel.

As shown in \eqref{def:tf}, the terminal time of CAV $i$ is dependent not only
on the terminal time of CAV $i-1$ and/or $k$, but also on the terminal speed
of CAV $i-1$ and/or $k$. Note that the terminal speed is unspecified and obtained from the energy
minimization problem \eqref{eq:decentral}. However, there is a number of ways to
specify the terminal speed.


\subsection{Alternative Energy Minimization Problem Formulations}

The effectiveness of the resequencing process may be affected by the way we
formulate the energy minimization problem. Next, by modifying
\eqref{eq:decentral}, we are going to explore several alternative problem
formulations and their impact on the resequencing efficiency.

\textbf{1) Modifying the terminal time of CAV \#1}, i.e., $t_1^m$: Due to
the recursive structure of the terminal times in \eqref{def:tf}, $t_{1}^{m}$
will generally affect all CAVs that follow CAV \#1. Recalling that there exists a
degree of freedom in the selection of $t_{1}^{m}$ which can be used to trade
off energy minimization and throughput maximization, we can modify
the energy minimization problem formulation for CAV \#1 by including the term
$\rho\cdot(t_{1}^{m}-t_{1}^{0})$ below to penalize longer travel times:
\begin{gather}
\min_{u_{1}}\frac{1}{2}\int_{t_{1}^{0}}^{t_{1}^{m}}u_{1}^{2}~dt+\rho(t_{1}%
^{m}-t_{1}^{0})\\
\text{subject to:~}\eqref{eq:model2},\eqref{speed_accel constraints},p_{1}(t_{1}^{m})=L,\text{given~}t_{1}^{0},v_{1}(t_{1}^{0}),p_{1}(t_{1}%
^{0}).\nonumber
\end{gather}
The coefficient $\rho$ allows trading off the throughput maximization and
energy minimization objectives. Note that the terminal time $t_{1}^{m}$ is now unspecified.
Alternatively, we can force CAV \#1 to reach the MZ as quickly as possible by
setting $t_{1}^{m}=t_{1}^{c}$, the lower bound for terminal times.

\textbf{2) Modifying the terminal speed of CAV} $\bm{i}$, i.e., $v_i^m$:
Due to the recursive terminal time structure in \eqref{def:tf}, the terminal
speed $v_{i}^{m}$ also impacts vehicles that follow $i$, hence, this affects
the traffic throughput. For example, a low terminal speed $v_{i-1}^{m}$ can
result in a long gap between CAV $i$ and $i-1$, which leads to a longer travel
time for $i$, thus reducing the traffic throughput. Therefore, we can modify
the energy minimization problem by including a quadratic deviation of $v_{i}^{m}$ from the maximum speed $v_{max}$ to penalize lower terminal
speeds, that is,
\begin{gather}
\min_{u_{i}}\frac{1}{2}\int_{t_{i}^{0}}^{t_{i}^{m}}u_{i}^{2}~dt+\frac{\sigma
}{2}(v_{i}^{m}-v_{max})^{2}\\
\text{subject to:~}%
\eqref{eq:model2},\eqref{speed_accel constraints}, t_i^m, p_{i}(t_{i}^{m})=L, \text{given~}t_{i}^{0},v_{i}(t_{i}^{0}),p_{i}(t_{i}^{0}).\nonumber
\end{gather}
The coefficient $\sigma$ allows trading off the throughput maximization and
energy minimization objectives.

Alternatively, we can directly set $v_{i}^{m}=v_{max}$. Note that CAV $i$ may
not be able to reach $v_{max}$. In that case, $v_{i}^{m}$ is set to the
maximal speed that CAV $i$ can reach given its initial conditions. Assuming
$v_{max}$ is reachable for CAV $i$, the energy minimization problem is
formulated as
\begin{gather}
\min_{u_{i}}\frac{1}{2}\int_{t_{i}^{0}}^{t_{i}^{m}}u_{i}^{2}~dt\\
\text{subject to:~}%
\eqref{eq:model2},\eqref{speed_accel constraints}, t_i^m, p_{i}(t_{i}^{m}), v_{i}^{m}=v_{max},t_{i}^{0},v_{i}(t_{i}^{0}),p_{i}(t_{i}^{0}).\nonumber
\end{gather}
Note that $v_{i}^{m}$ is specified in this formulation.


\subsection{Case Study for Dynamic Resequencing}

The effectiveness of the resequencing process in terms of maximizing the
traffic throughput is validated through simulation in MATLAB considering 20
CAVs crossing an urban intersection. The intersection is asymmetric by setting
the lengths of the CZ segments to $L_{E2W}=L_{W2E}=400$m and $L_{N2S}%
=L_{S2N}=300$m, respectively. The width of the merging zone is $S=30$m. The
vehicle arrival process is assumed to be given by a Poisson process with the
same rate $\lambda=0.4$ (veh/s) for each CZ
segment. The initial speeds are assumed to be given by a uniform distribution
defined over $[8,12]$ m/s. The maximum speed and maximum acceleration are
$v_{max}=16$ m/s and $u_{max}=2$ m/s$^{2}$, while the minimum speed and
maximum deceleration (minimum acceleration) are set to $v_{min}=4$ m/s and
$u_{min}=-5$ m/s$^{2}$.

We consider 10 different alternative energy minimization problem formulations
for comparison ([R] indicates a case with resequencing, and [NR] without resequencing):

\begin{itemize}
\item[(1)] [NR] CAV \#1 cruises and reaches MZ at $t_{1}^{m}
= t_{1}^{0} + \frac{L}{v_{i}^{0}}$;

\item[(2)] [NR] CAV \#1 is penalized for longer travel time by including the
term $\rho(t_{1}^{m} - t_{1}^{0})$ in the cost functional, where $\rho= 5$;

\item[(3)] [NR] CAV \#1 is forced to reach MZ at $t_{1}^{m} = t_{1}^{c}$;

\item[(4)] [R] CAV \#1 cruises and reaches MZ at $t_{1}^{m}
= t_{1}^{0} + \frac{L}{v_{i}^{0}}$;

\item[(5)] [R] CAV \#1 is penalized for longer travel time by including the term
$\rho(t_{1}^{m} - t_{1}^{0})$ in the cost functional, where $\rho= 5$;

\item[(6)] [R] CAV \#1 is forced to reach MZ at $t_{1}^{m} = t_{1}^{c}$;

\item[(7)] [R] CAVs are penalized from deviating $v_{max}$ at $t_{i}^{m}$ by including the term $\frac{\sigma}{2}(v_{i}(t) - v_{max})^{2}$ in
the cost functional, where $\sigma= 0.1$; 

\item[(8)] [R] similar to case (7), except that $\sigma= 1$; 

\item[(9)] [R] similar to case (7), except that $\sigma= 10$; 

\item[(10)] [R] CAVs are forced to reach $v_{max}$ at $t_{i}^{m}$.
\end{itemize}

\begin{figure}[tbh]
\centering
\includegraphics[width=1\columnwidth]{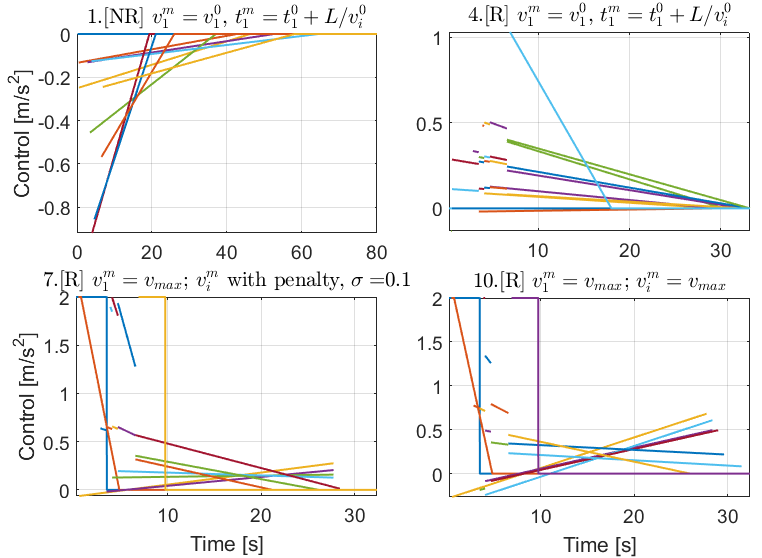} \caption{Optimal
control profiles of the first 10 CAVs under part of the problem formulation cases.}%
\label{uit}%
\end{figure}

\begin{figure}[tbh]
\centering
\includegraphics[width=1\columnwidth]{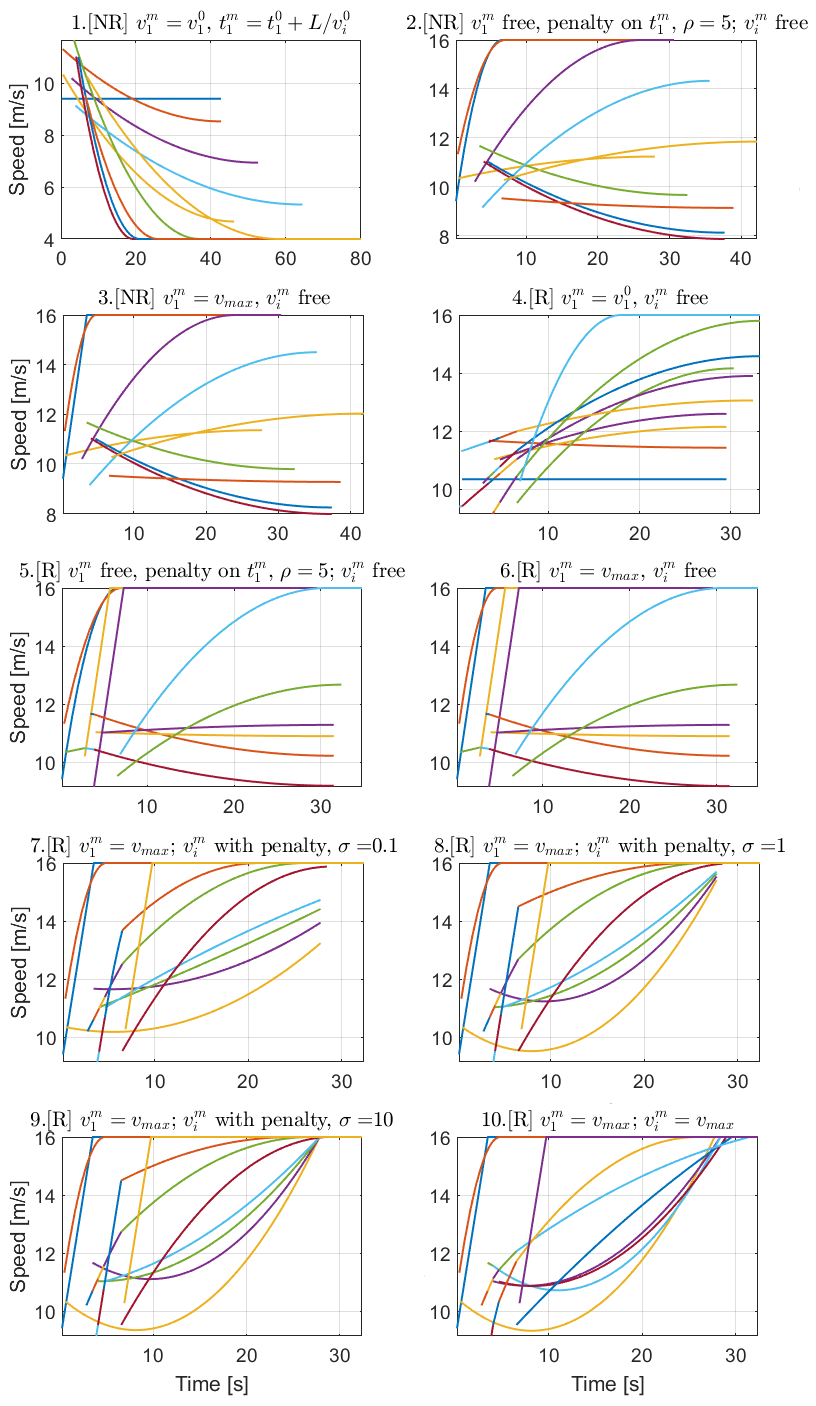} \caption{Optimal speed
trajectories of the first 10 CAVs under different problem formulation cases.}%
\label{vit}
\end{figure}

The optimal control and speed trajectories of the first 10 CAVs under different problem formulation cases are shown in
Fig. \ref{uit} and \ref{vit} respectively. Within each trajectory, the change
of color indicates an occurrence of a resequencing process. In Fig. \ref{uit},
observe that there may exist a discontinuity within a control trajectory when
the resequencing process takes place since resequencing may lead to an updated
optimal trajectory. Note that the speed and control constraints \eqref{speed_accel constraints} are satisfied throughout the trajectories. 

\begin{figure}[tbh]
\centering
\includegraphics[width=0.9\columnwidth]{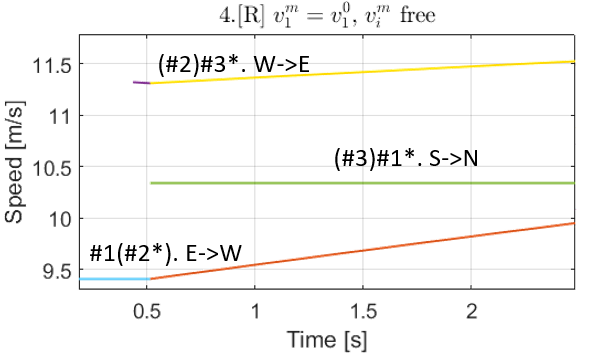} \caption{Illustration of the resequencing process.}%
\label{resequencing_example}%
\end{figure}

To illustrate the resequencing process, part of the
speed trajectories for the first 3 CAVs under case 4 are shown in Fig. \ref{resequencing_example}, where CAV \#1 is
cruising in an energy-optimal way, i.e., $t_{1}^{m}=t_{1}^{0}+\frac{L}%
{v_{i}^{0}}$ and the crossing sequence is re-evaluated whenever a CAV enters
the CZ. Observe that when CAV \#3 arrives at the CZ, it is rescheduled to
\#1$_{3}$ (previously indexed as \#3), and CAV \#1 and \#2 are rescheduled to
\#2$_{1}$ and \#3$_{2}$. 
Note that both CAV \#2 and \#1 are traveling on the longer CZ segments, while CAV \#3 is traveling on the shorter CZ segment. Intuitively, since CAV \#3 enters the CZ right after CAV \#2
($t_{2}^{0}=0.43$s, and $t_{3}^{0}=0.51$s), it is natural to let CAV \#3 cross the MZ first as it is closer to the MZ. Without
the resequencing process (case 1), CAV \#3 can only enter the MZ when \#2
leaves the MZ, which makes the total gap $t_{3}^{m}-t_{1}^{m}=3.52$s; with the
resequencing process (case 4), CAV \#3 becomes \#1$_{3}$, and the total gap
reduces to $t_{3_{2}}^{m}-t_{1_{3}}^{m}=2.9$s, hence, improving the traffic throughput.

Under case 1, where CAV \#1 is assumed to cruise at its initial speed in terms
of minimizing energy consumption and no resequencing is considered, CAV\#3
results in a low terminal speed $v_{3}^{m}=4.67$m/s. Under case 4 where the
resequencing process is included, CAV \#3 is rescheduled to \#1$_{3}$ and
assumed to cruise at its initial speed. Therefore, the terminal times for CAV
\#2$_{1}$ and \#3$_{2}$ are updated based on the recursive terminal time
structure and observe that $t_{3_{2}}^{m}<t_{3}^{m}$. This forces CAV
\#3$_{2}$ to accelerate, which leads to a higher terminal speed $v_{3_{2}}%
^{m}=13.1$m/s and further minimizes the gap.

\begin{customrem}{1}
This case study assumes a vehicle arrival rate near the saturation level
(further discussed in Sec. \ref{bigo}), which indicates that the gaps between
CAV arrivals at the CZ, i.e., $t_{i}^{0}-t_{i-1}^{0}$, are relatively small. When the gap between
CAV arrivals is smaller than the gap between
terminal times, i.e., $t_{i}^{m}-t_{i-1}^{m}$, CAV $i$ is naturally forced to
slow down as the terminal speed is undefined in \eqref{eq:decentral} and
results in lower terminal speed. Conversely, when the gap between CAV arrivals
is larger than the gap between terminal times, CAV $i$ may need to accelerate
which leads to higher terminal speed. \label{remark:1}
\vspace{-4mm}
\end{customrem}

Since the resequencing process aims at finding the optimal crossing sequence
which maximizes the traffic throughput, the cases with the resequencing (Fig. \ref{vit}(4-10)) outperform those without resequencing (Fig.
\ref{vit}(1-3)). This can be seen by comparing the total travel time
among different cases. In addition, due to the recursive structure of the
terminal times, there exists a propagation effect of the terminal speeds: a
lower terminal speed of CAV $i-1$ may lead to higher terminal time for CAV
$i$, which further lowers the terminal speed of $i$, as shown in Fig.
\ref{vit}(1). With the resequencing process, CAV $i$ may be rescheduled to an
earlier position $j_{i}<i$ in the queue. Therefore, $t_{j_{i}}^{m}$
$<t_{i}^{m}$, which leads to higher terminal speed $v_{j_{i}}^{m}$. Even
though CAV $j$ (now indexed as $(j+1)_{j}$) is affected by the resequencing
process, the increase in $t_{(j+1)_{j}}^{m}$ is minimal due to the higher
$v_{j_{i}}^{m}$. Thus, the gap decreases and the traffic throughput
improves compared to the cases without resequencing.

In cases 6 to 10 (Fig. \ref{vit}(6-10)), we are increasing the weight forcing
the terminal speeds of CAVs to reach $v_{max}$. Note that the travel times in these cases are similar due to the fact that resequencing  results in lower terminal times, which naturally
leads to higher terminal speeds even without forcing a CAV to reach $v_{max}$.

Observe that without the resequencing process (Fig. \ref{vit}(1-3)), changing
$t_{1}^{m}$ alone can affect the traffic throughput. In Fig. \ref{vit}(1), CAV
\#1 is assumed to be cruising at its initial speed. Due to the propagation
effect of the terminal speeds, the following CAVs end up with lower terminal
speeds, which decreases the total travel time. In Fig. \ref{vit}(2-3), as we
are forcing CAV \#1 to reach $v_{max}$ when it arrives at the MZ, the terminal
time $t_{1}^{m}$ decreases, hence, $t_{i}^{m}$, $i>1$, determined by the
recursive terminal time structure, also decreases. Thus, the following
vehicles result in higher terminal speeds, which reduces the total travel time
by a large margin. The benefit obtained from varying $t_{1}^{m}$ diminishes in
the cases with resequencing (Fig. \ref{vit}(4-6)). 


\subsection{Performance Metrics}
\label{performance}
\begin{figure}[tbh]
\centering
\includegraphics[width=1 \columnwidth]{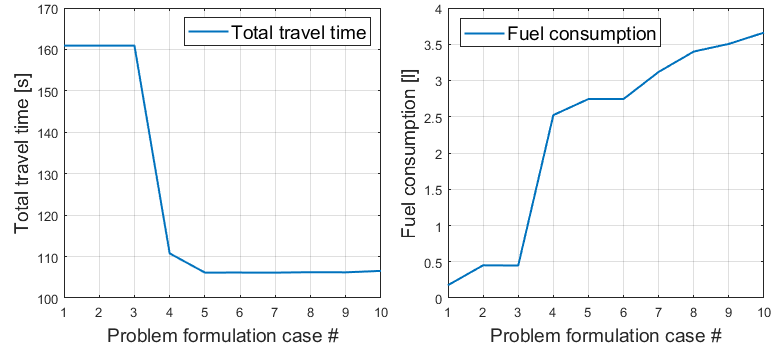}
\caption{Travel time (left) and fuel consumption (right) under alternative
problem formulations given $\lambda= 0.4$ (veh/s).}%
\label{metric}%
\end{figure}

To quantify the effectiveness of the resequencing process, we compare the
performance metrics, i.e., energy consumption and throughput under different
cases. To measure the throughput, we use $t_{N(t)}^{m}$, the time by which all
$N(t)$ vehicles exit the CZ. To measure the energy consumption, we use the
polynomial metamodel in \cite{Kamal2013a} that yields vehicle fuel consumption
as a function of speed and acceleration. We consider 100 CAVs
crossing one intersection given a vehicle arrival rate of $\lambda=0.4$ (veh/s).
The performance metrics are shown in Fig.
\ref{metric}.  Observe that with the
resequencing process (starting with case 4), the travel time is improved by
approximately 34\% compared to the cases without resequencing. This is
consistent with the observations discussed in the case study and shows the
efficiency of the resequencing process in terms of traffic throughput maximization.

In contrast to what we have observed in Fig. \ref{vit}(1-3), cases 1, 2, and 3 achieve almost the same
travel time in Fig. \ref{metric}. This leads to the conclusion that how we specify $v_{1}^{m}$ does not have any effect when traffic flows reach steady state.

In Fig. \ref{metric}, observe that the resequencing process leads to an increase in energy, counteracting the throughput benefits. This shows the tradeoff between energy minimization and throughput maximization. Unlike the travel time curve where cases 4-10 have minimal difference in improving travel time, as we are increasingly forcing terminal speeds to reach $v_{max}$ (from case 4 to 10), more fuel is consumed. As a whole, cases 4 and 5 achieve better performance compared to others.

To further investigate the tradeoff between the throughput maximization and
energy minimization objectives, we explore the two performance metrics over
cases 4-10 when resequencing is applied at different traffic intensities, as
summarized in Table \ref{tab:volume}. Observe that as the traffic intensity
decreases, the average travel time is improved, while more fuel is expended.
Also observe that when the traffic is light, e.g., $\lambda=0.1$ (veh/s), the average
travel times do not significantly vary over different problem formulations.
Due to the light traffic, the recursive terminal time structure is interrupted
by the critical time $t_{i}^{c}$. Generally, lower travel time corresponds to
more fuel consumption, which is consistent with the expected tradeoff between
energy minimization and throughput maximization.

\begin{table}[ht]
\caption{Performance under different traffic intensities}%
\label{tab:volume}%
\centering
\begin{tabular}
{|p{5mm}<{\centering}|p{5mm}<{\centering}|p{5mm}<{\centering}|p{5mm}<{\centering}|p{5mm}<{\centering}|p{5mm}<{\centering}|p{5mm}<{\centering}|p{5mm}<{\centering}|p{5mm}<{\centering}|}\hline 
\multirow{2}*{} &  & 4 & 5 & 6 & 7 & 8 & 9 & 10\\\hline
\multirow{2}*{$\lambda$=0.4} & time & 32.44 & 28.11 & 28.1 & 28 & 28.1 &
28.16 & 29.4\\\cline{2-9}
& fuel & 1.55 & 2.09 & 2.09 & 2.19 & 2.26 & 2.28 & 2.25\\\hline
\multirow{2}*{$\lambda$=0.3} & time & 28.53 & 26.98 & 26.99 & 27.04 & 27.38 &
27.42 & 28.17\\\cline{2-9}
& fuel & 2.02 & 2.2 & 2.2 & 2.25 & 2.29 & 2.3 & 2.29\\\hline
\multirow{2}*{$\lambda$=0.2} & time & 27.19 & 26.43 & 26.40 & 26.46 & 26.73 &
26.72 & 26.91\\\cline{2-9}
& fuel & 2.17 & 2.25 & 2.26 & 2.28 & 2.31 & 2.32 & 2.31\\\hline
\multirow{2}*{$\lambda$=0.1} & time & 26.17 & 25.98 & 26.04 & 26.27 & 26.30 &
26.29 & 26.33\\\cline{2-9}
& fuel & 2.26 & 2.28 & 2.28 & 2.29 & 2.33 & 2.34 & 2.34\\\hline
\multicolumn{9}{l}{$\lambda$: arrival rate in veh/s; time in second;
fuel in liter}
\end{tabular}
\vspace{-3mm}
\end{table}

\begin{customrem}{2}
The terminal times are recursively computed based on the lateral and rear-end
collision avoidance constraints. These safety constraints are conservative in
the sense that only one vehicle is allowed inside the MZ if CAV $i-1\in
\mathcal{C}_{i}(t)$. However, the traffic throughput can always be improved by
subdividing the MZ into smaller single-vehicle areas and establishing less
conservative safety constraints. \label{remark:2}
\end{customrem}


\section{Computational Complexity Analysis for Resequencing}

\label{bigo} Since the coordinator needs to re-evaluate the crossing sequence
every time a new CAV arrives at the CZ, the complexity of the resequencing
process (see Algorithm 1) may be significant when the traffic is heavy. A key
observation is that CAV $i$ can obviously not overtake its preceding CAV $k$,
which therefore, guarantees an upper bound in the resequencing computational
complexity involved. Since the key to the resequencing process lies in
inserting CAV $i$ into different positions of the queue after $k$, the
computational complexity can be represented by the number of swaps $f(i,i-1)$
in addition to the computation without resequencing.

In what follows, we carry out first a worst case analysis. This corresponds to CAV
$i$ entering the CZ when there is no preceding vehicle $k$ traveling on the
same lane, while all other CZ road segments operate near capacity. Assuming
four CZ segments within an intersection, their lengths are denoted by $L_{r}$,
$r\in\{1,2,3,4\}$. Vehicle arrivals are assumed to be distributed according to
Poisson processes with rates $\lambda_{r}$, $r\in\{1,2,3,4\}$. Letting the
average CAV length be $l_{v}$, the capacity for each CZ segment $C_{r}$ is
given by $
C_{r}=\frac{L_{r}}{l_{v}+\delta},\text{ \ \ }r\in\{1,2,3,4\}$. Assuming CAV $i$ enters the first CZ segment, i.e., $r=1$, the computational complexity measured using the number of swaps
for CAV $i$, denoted as $N$, under the worst case is $N^{1}=\frac{L_{2}+L_{3}+L_{4}}{l_{v}+\delta} + 1$. Taking the vehicle arrivals on other CZ segments into consideration, the worst
case of the computational complexity for the whole intersection is $N=\text{max}\{\frac{L_{2}+L_{3}+L_{4}}{l_{v}+\delta},\frac{L_{1}%
+L_{3}+L_{4}}{l_{v}+\delta},\frac{L_{1}+L_{2}+L_{4}}{l_{v}+\delta},\frac{L_{1}+L_{2}+L_{3}}{l_{v}+\delta
}\} + 1$.
This represents the \emph{upper bound} of the computational complexity
associated with the resequencing process. The best case occurs when $k=i-1$,
which indicates no necessity to resequence. Hence, the \emph{lower bound} is
$N=1$.

The \textit{saturation flow rate} is an important concept associated with the
stability of the intersection viewed as a queueing system. When the
intersection is saturated, the number of vehicles present exceeds its capacity
and congestion occurs. In this case, it is not possible to apply any control
other than traffic signaling. Thus, it is important to derive the
\emph{expected} computational complexity when the traffic flow is stable. The
saturation flow rate is defined as the headway (in time units) between
vehicles moving at steady state. Viewed as a queueing system, the intersection
is an M/G/1 queue, where the MZ is the server and the vehicles are the
customers in the queue. The condition for this M/G/1 queueing system to be
stable is $\sum_{r\in\{1,2,3,4\}}\lambda_{r}<\mu$, where $\lambda_{r}$ is the arrival rate on $r$th road
segment, and $\mu$ is the service rate of the MZ. Based on the recursive
structure of terminal times in \eqref{def:tf}, vehicles traveling on
opposite roads will not generate any collision inside the MZ, hence, they
are allowed to cross the MZ at the same time. It follows that we only need
$\sum_{r\in\{1,2,3,4\}}\lambda_{r}<2\mu$ as a condition for stability.

\textit{Expected computational complexity} $E[N]$: to compute $E[N]$, we first consider the expected interarrival time between
CAVs $k$ and $i$. Assuming that CAV $i$ enters the first CZ segment, i.e., $r=1$, the
expected interarrival time is $E[\Delta t]=\frac{1}{\lambda_{1}}$. Over the
interarrival time $\Delta t$, the expected number of arrivals on the other three CZ
segments are given by $E[\Delta t]\cdot(\lambda_{2}+\lambda_{3}+\lambda_{4})$.
Therefore, for vehicles coming from the first CZ segment, we have 
\[E[N^{1}]=\frac{\lambda_{2}+\lambda_{3}+\lambda_{4}}{\lambda_{1}}+1.
\]
Similarly, for the other three CZ segments, we have $E[N^{2}]=\frac
{\lambda_{1}+\lambda_{3}+\lambda_{4}}{\lambda_{2}}+1$, $E[N^{3}]=\frac
{\lambda_{1}+\lambda_{2}+\lambda_{4}}{\lambda_{3}}+1$, $E[N^{4}]=\frac
{\lambda_{1}+\lambda_{2}+\lambda_{3}}{\lambda_{4}}+1$. Therefore,
\begin{gather}
\text{E}[N]  = \frac{\lambda_1  \text{E}[N^1] + \lambda_2  \text{E}[N^2] + \lambda_3  \text{E}[N^3]+ \lambda_4  \text{E}[N^4]}{(\lambda_1 + \lambda_2 + \lambda_3 + \lambda_4)} = 4  \nonumber 
\end{gather}
regardless of the arrival rates. Thus, the expected computational complexity
$E[N]=4$ happens to be the number of CZ segments. In fact, this result can be
generalized to an intersection with $M$ lanes: for vehicles coming from the
$p$th CZ segment, the expected computational complexity $E[N^{p}]$ can be shown to be
\[
E[N^{p}]=\frac{1}{\lambda_{p}}\sum_{r\in\{1,\dots,M\},r\neq p}\lambda_{r}+1,
\]
and for the whole intersection, we have
\begin{gather}
 \text{E}[N] = \frac{\lambda_1 \cdot \text{E}[N^1] + \dots + \lambda_M \cdot \text{E}[N^M]}{\sum_{r \in \{1, \dots, M\}} \lambda_r}  = M  
\end{gather}
This indicates that the expected computational complexity is always determined
by the number of lanes associated with the intersection.

The expected computational complexity is validated through simulation in
MATLAB considering 100 CAVs crossing an urban intersection, with exactly the same simulation settings as in Sec. \ref{performance} with $M=4$ lanes. The average service time is roughly estimated as 1.25s and the expected service rate is $\mu=0.8$. Therefore, the stability condition can be determined as $\lambda_{1}+\lambda_{2}+\lambda_{3}+\lambda_{4}<1.6$.
The energy minimization problem is formulated as in case 5 in
Section III.D, which penalizes longer travel times for CAV \#1.

\begin{figure}[tbh]
\centering
\includegraphics[width=1\columnwidth]{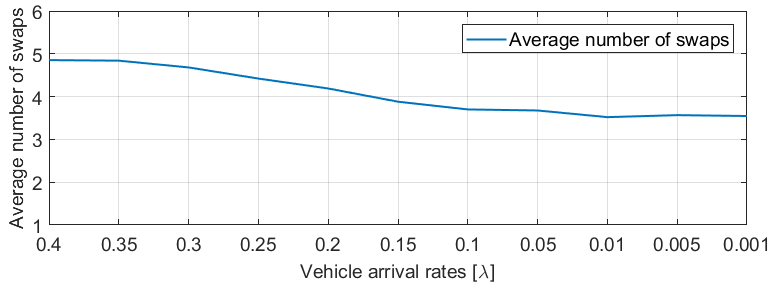} \caption{Expected computational complexity of resequencing process over decreasing traffic intensity.}%
\label{lambda}%
\end{figure}


The simulation results are shown in Fig. \ref{lambda}, where the computational complexity, measured using the number of swaps, is averaged over 10 simulations. We assume $\lambda_{1}=\lambda_{2}=\lambda_{3}=\lambda_{4} = \lambda$, where $\lambda < 0.4$. Over different arrival rates, the computational complexity in performing dynamic resequencing, is approximately 4, as expected. 
The actual value of $E[N]$, however, may be lower since a
resequencing step affects subsequent resequencing steps by altering the
vehicle arrival process distribution.

\section{CONCLUSIONS AND FUTURE WORK}

Earlier work \cite{ZhangMalikopoulosCassandras2016}, \cite{Malikopoulos2016} and \cite{Zhang2016} has established a decentralized optimal control framework for optimally controlling CAVs crossing a signal-free urban
intersection while following a strict FIFO queueing order. In this paper, we
have extended the solution of this problem to account for asymmetric
intersections by relaxing the FIFO constraint and introducing a dynamic
resequencing process so as to maximize the traffic throughput. The dynamic
resequencing has been shown to be computationally very efficient. It is also
shown to reduce the travel time at the cost of additional fuel consumption.
This tradeoff has been illustrated  through simulation examples.

Ongoing research is considering turns (see \cite{Zhang2017}) and lane changing
in the intersection with a diverse set of CAVs and exploring the mixed
scenario where both CAVs and human-driven vehicles travel on the roads (see \cite{Zhang2018}). Future
research should also investigate the multi-intersection scenario and how the
coupling between multiple intersections would affect the throughput
maximization and energy minimization problems.

\vspace{-1mm}
\bibliographystyle{IEEETran}
\bibliography{CDC2018}

\begin{thebibliography}{10}
\providecommand{\url}[1]{#1}
\csname url@rmstyle\endcsname
\providecommand{\newblock}{\relax}
\providecommand{\bibinfo}[2]{#2}
\providecommand\BIBentrySTDinterwordspacing{\spaceskip=0pt\relax}
\providecommand\BIBentryALTinterwordstretchfactor{4}
\providecommand\BIBentryALTinterwordspacing{\spaceskip=\fontdimen2\font plus
\BIBentryALTinterwordstretchfactor\fontdimen3\font minus
  \fontdimen4\font\relax}
\providecommand\BIBforeignlanguage[2]{{%
\expandafter\ifx\csname l@#1\endcsname\relax
\typeout{** WARNING: IEEEtran.bst: No hyphenation pattern has been}%
\typeout{** loaded for the language `#1'. Using the pattern for}%
\typeout{** the default language instead.}%
\else
\language=\csname l@#1\endcsname
\fi
#2}}

\bibitem{Fleck2015}
J.~L. Fleck, C.~G. Cassandras, and Y.~Geng, ``Adaptive quasi-dynamic traffic
  light control,'' \emph{IEEE Transactions on Control Systems Technology},
  2015, DOI: 10.1109/TCST.2015.2468181, to appear.

\bibitem{Levine1966}
W.~Levine and M.~Athans, ``{On the optimal error regulation of a string of
  moving vehicles},'' \emph{IEEE Transactions on Automatic Control}, vol.~11,
  no.~3, pp. 355--361, 1966.

\bibitem{Dresner2004}
K.~Dresner and P.~Stone, ``{Multiagent traffic management: a reservation-based
  intersection control mechanism},'' in \emph{Proceedings of the Third
  International Joint Conference on Autonomous Agents and Multiagents Systems},
  2004, pp. 530--537.

\bibitem{Dresner2008}
------, ``{A Multiagent Approach to Autonomous Intersection Management},''
  \emph{Journal of Artificial Intelligence Research}, vol.~31, pp. 591--653,
  2008.

\bibitem{DeLaFortelle2010}
A.~{de La Fortelle}, ``{Analysis of reservation algorithms for cooperative
  planning at intersections},'' \emph{13th International IEEE Conference on
  Intelligent Transportation Systems}, pp. 445--449, Sept. 2010.

\bibitem{Huang2012}
S.~Huang, A.~Sadek, and Y.~Zhao, ``{Assessing the Mobility and Environmental
  Benefits of Reservation-Based Intelligent Intersections Using an Integrated
  Simulator},'' \emph{IEEE Transactions on Intelligent Transportation Systems},
  vol.~13, no.~3, pp. 1201,1214, 2012.

\bibitem{Li2006}
L.~Li and F.-Y. Wang, ``{Cooperative Driving at Blind Crossings Using
  Intervehicle Communication},'' \emph{IEEE Transactions in Vehicular
  Technology}, vol.~55, no.~6, pp. 1712,1724, 2006.

\bibitem{Yan2009}
F.~Yan, M.~Dridi, and A.~{El Moudni}, ``{Autonomous vehicle sequencing
  algorithm at isolated intersections},'' \emph{2009 12th International IEEE
  Conference on Intelligent Transportation Systems}, pp. 1--6, 2009.

\bibitem{Zhu2015}
F.~Zhu and S.~V. Ukkusuri, ``{A linear programming formulation for autonomous
  intersection control within a dynamic traffic assignment and connected
  vehicle environment},'' \emph{Transportation Research Part C: Emerging
  Technologies}, 2015.

\bibitem{Miculescu2014}
D.~Miculescu and S.~Karaman, ``{Polling-Systems-Based Control of
  High-Performance Provably-Safe Autonomous Intersections},'' in \emph{53rd
  IEEE Conference on Decision and Control}, 2014.

\bibitem{gilbert1976vehicle}
E.~G. Gilbert, ``Vehicle cruise: Improved fuel economy by periodic control,''
  \emph{Automatica}, vol.~12, no.~2, pp. 159--166, 1976.

\bibitem{hooker1988optimal}
J.~Hooker, ``Optimal driving for single-vehicle fuel economy,''
  \emph{Transportation Research Part A: General}, vol.~22, no.~3, pp. 183--201,
  1988.

\bibitem{hellstrom2010design}
E.~Hellstr{\"o}m, J.~{\AA}slund, and L.~Nielsen, ``Design of an efficient
  algorithm for fuel-optimal look-ahead control,'' \emph{Control Engineering
  Practice}, vol.~18, no.~11, pp. 1318--1327, 2010.

\bibitem{li2012minimum}
S.~E. Li, H.~Peng, K.~Li, and J.~Wang, ``Minimum fuel control strategy in
  automated car-following scenarios,'' \emph{IEEE Transactions on Vehicular
  Technology}, vol.~61, no.~3, pp. 998--1007, 2012.

\bibitem{Rios-Torres}
J.~Rios-Torres and A.~A. Malikopoulos, ``{A Survey on the Coordination of
  Connected and Automated Vehicles at Intersections and Merging at Highway
  On-Ramps},'' \emph{IEEE Transactions on Intelligent Transportation Systems},
  2016 (forthcoming).

\bibitem{ZhangMalikopoulosCassandras2016}
Y.~Zhang, A.~A. Malikopoulos, and C.~G. Cassandras, ``Optimal control and
  coordination of connected and automated vehicles at urban traffic
  intersections,'' in \emph{Proceedings of the American Control Conference},
  2016, pp. 6227--6232.

\bibitem{Malikopoulos2016}
A.~A. Malikopoulos, C.~G. Cassandras, and Y.~Zhang, ``A decentralized
  energy-optimal control framework for connected automated vehicles at
  signal-free intersections,'' \emph{Automatica}, 2017, (to appear).

\bibitem{Zohdy2012}
I.~H. Zohdy, R.~K. Kamalanathsharma, and H.~Rakha, ``{Intersection management
  for autonomous vehicles using iCACC},'' \emph{2012 15th International IEEE
  Conference on Intelligent Transportation Systems}, pp. 1109--1114, 2012.

\bibitem{Lee2012}
J.~Lee and B.~Park, ``{Development and Evaluation of a Cooperative Vehicle
  Intersection Control Algorithm Under the Connected Vehicles Environment},''
  \emph{IEEE Transactions on Intelligent Transportation Systems}, vol.~13,
  no.~1, pp. 81--90, 2012.

\bibitem{Zhang2016}
Y.~Zhang, C.~G. Cassandras, and A.~A. Malikopoulos, ``Optimal control of
  connected automated vehicles at urban traffic intersections: A feasibility
  enforcement analysis,'' in \emph{Proceedings of the 2017 American Control
  Conference}, 2017, pp. 3548--3553.

\bibitem{Kamal2013a}
M.~Kamal, M.~Mukai, J.~Murata, and T.~Kawabe, ``{Model Predictive Control of
  Vehicles on Urban Roads for Improved Fuel Economy},'' \emph{IEEE Transactions
  on Control Systems Technology}, vol.~21, no.~3, pp. 831--841, 2013.

\bibitem{Zhang2017}
Y.~Zhang, A.~A. Malikopoulos, and C.~G. Cassandras, ``Decentralized optimal
  control for connected automated vehicles at intersections including left and
  right turns,'' in \emph{56th IEEE Conference on Decision and Control}, 2017,
  pp. 4428--4433.

\bibitem{Zhang2018}
Y.~Zhang and C.~G. Cassandras, ``The penetration effect of connected automated
  vehicles in urban traffic: an energy impact study,'' 2018, arXiv: 1803.05577.

\end{thebibliography}

\end{document}